\documentclass{interact}
\makeatletter
\def\ps@pprintTitle{%
	\let\@oddhead\@empty
	\let\@evenhead\@empty
	\def\@oddfoot{\centerline{\thepage}}%
	\let\@evenfoot\@oddfoot}
\makeatother

\usepackage{amsfonts}
\usepackage{enumerate} 
\usepackage{amsthm}
\usepackage{amsmath}
\usepackage{mathtools}
\usepackage{amscd}
\usepackage[active]{srcltx} 
\usepackage{latexsym}
\usepackage{amssymb}
\usepackage{enumitem}

\newtheorem{thm}{Theorem}[section]
\newtheorem{prop}[thm]{Proposition}
\newtheorem{cor}[thm]{Corollary}
\newtheorem*{cor*}{Corollary}
\newtheorem{lema}[thm]{Lemma}
\newtheorem*{lema*}{Lemma}

\numberwithin{equation}{section}
\theoremstyle{definition}
\newtheorem*{Def}{Definition}
\newtheorem*{Dem}{Proof}
\newenvironment{dem}{\vspace{1ex}\noindent{\it Proof.}\hspace{0.5em}}
{\hfill\qed\vspace{1ex}}

\newtheorem*{obs*}{Remark}
\newtheorem*{thm*}{Theorem}
\newtheorem*{prop*}{Proposition}

\newtheoremstyle{dotless}{}{}{}{}{}{}{ }{}
\theoremstyle{dotless}

\newcommand{\PI}[2]{\left\langle \,#1 , #2\, \right\rangle}
\newcommand{\PIW}[2]{\left\langle \,#1 , #2\, \right\rangle_W}

\newcommand{\NWW}[1]{\Vert #1 \Vert_{W}}

\newcommand{\ra}{\rightarrow}

\newcommand{\WN}{W_{/ N(W_{/ \St})}}
\newcommand{\TT}{N(W_{/ \St})}
\newcommand{\TTA}{N(W_{/ R(A)})}

\def\minus{\stackrel{-}{\leq}}
\def\star{\stackrel{*}{\leq}}
\def\starW{\stackrel{*_W}{\leq}}

\def\lminus{{\;{\vphantom{\leq}}_{-}\leq\,}}
\def\rminus{{\;\leq{\vphantom{\leq}}_{-}\;}}

\def\lW{{\;{\vphantom{\leq}}_{*_W}\leq\,}}

\def\lstar{{\;{\vphantom{\leq}}_{*}\leq\,}}

\def\rW{{\;\leq{\vphantom{\leq}}_{*_W}\;}}
\def\rstar{\;{\leq{\vphantom{\leq}}_{*}\;}}

\newcommand{\St}{\mathcal{S}}

\newcommand{\HH}{\mathcal{H}}

\newcommand{\M}{\mathcal{M}}
\newcommand{\N}{\mathcal{N}}

\newcommand{\Q}{\mathcal{Q}}

\newcommand{\mc}[1]{\mathcal{#1}}

\newcommand{\ol}{\overline}

\begin{document}

\title{Shorted operators and minus order}
\author{
\name{M. Contino\textsuperscript{a,b}
J. I. Giribet\textsuperscript{a,b}
A. Maestripieri\textsuperscript{a,b}
}
\affil{
\textsuperscript {a} Instituto Argentino de Matem\'atica "Alberto Calder\'on" - CONICET. Saavedra 15, piso 3 (1083) Ciudad Aut\'onoma de Buenos Aires, Argentina.\\ 
\textsuperscript {b} Facultad de Ingenier\'ia - Universidad de Buenos Aires. Paseo Col\'on 850 (1063), Ciudad Aut\'onoma de Buenos Aires, Argentina.
}
}

\maketitle

\begin{abstract}

Let $\HH$ be a Hilbert space, $L(\HH)$ the algebra of bounded linear operators on $\HH$ and $W \in L(\HH)$ a positive operator. Given a closed subspace $\St$ of $\HH$, we characterize the shorted operator $W_{/ \St}$  of $W$ to $\St$ as the maximum and as the infimum of certain sets, for the {\it minus order} $\minus$. 
Also, given $A \in L(\HH)$ with closed range, we study the following operator approximation problem considering the minus order: 
$$
min_{\minus} \  \{(AX-I)^*W(AX-I) :  X \in L(\HH) ,\label{eqa2} \mbox{ subject to } N(A^*W)\subseteq N(X) \}.
$$

We show that, under certain conditions, the shorted operator of $W_{/R(A)}$ is the minimum of this problem and we characterize the set of solutions. 

\end{abstract}


\section{Introduction}

The minus order was introduced by Hartwig \cite{Har80} and independently by Nambooripad \cite{Nam80}, in both cases on semigroups, with the idea of generalizing some classical partial orders.  In \cite{AntCorSto06}, Antezana et al., extended the notion of the minus order to bounded linear operators acting on an infinite-dimensional Hilbert space. In this case, if $\HH$ is a Hilbert space and $L(\HH)$ is the algebra of bounded linear operators acting on $\HH,$ for $A, B \in L(\HH),$ $A \minus B$ (where the symbol $\minus$ stands for the minus order of operators) if the \emph{Dixmier} angle between $\ol{R(A)}$ and $\ol{R(B-A)}$ and the Dixmier angle between $\ol{R(A^*)}$ and $\ol{R(B^*-A^*)}$ are less than 1, where $R(T)$ stands for the range of the operator $T.$ Independently, in \cite{Sem10} $\check{\textrm{S}}$emrl gave another characterization of the minus order, he showed that $A \minus B$  if and only if there are bounded (oblique) projections, i.e. idempotents, $P$ and $Q$ such that $A=PB$ and $A^* =QB^*.$ 
In \cite{Minus}, a  new characterization of the minus order for operators acting on Hilbert spaces was given in terms of the so called range additivity property. Namely, it was proved that $A\minus B$ is equivalent to the range of $B$ being the direct sum of ranges of $A$ and $B-A$ and the range of $B^*$ being the direct sum of ranges of $A^*$ and $B^*-A^*$, which generalizes previous results presented in \cite{Sem10}. This plays an equivalent role to the rank additivity characterization when $A$ and $B$ are matrices \cite{Har80,Mitramatrixpartial}.

In \cite{AntCorSto06} the notion of shorted operator appears in relation with the minus order.
Given  a closed subspace $\St$ of $\HH$ and $W \in L(\HH)$ a positive operator, in 1947,  Krein \cite{Krein}, proved the existence of a maximum (with respect to the  order induced by the cone of positive operators) of the set $$\M(W,\St)=\{ X \in L(\HH): \ 0\leq X\leq W \mbox{ and } R(X)\subseteq \St^{\perp}\}.$$ Krein used this extremal operator in his theory of extension of symmetric operators. Years later, Anderson and Trapp \cite{Shorted2} studying the same problem, called this maximum the {\it shorted operator of $W$ to $\St$} (in the following denoted $W_{ / \St }$) and showed interesting properties of this operator and its connections with electrical circuit theory. The shorted operator have shown to be useful in many applications \cite{AntCorSto06}, \cite{Zhang}.

The pair $(W,\St)$ is said to be \emph{compatible} if there exists a bounded linear (not necessarily selfadjoint) projection $Q$ onto $\St$ such that $WQ$ is selfadjoint. Thus, if
$$\mc{P}(W,\St)=\{ Q \in L(\HH): Q^2=Q, \ R(Q)=\St, \ WQ=Q^*W \},$$ then $(W,\St)$ is compatible if and only if $\mc{P}(W,\St)$ is not empty. In \cite{CMSSzeged}, it was shown that there exists a strong relationship between compatibility, the projections of $\mc{P}(W,\St)$ and the shorted operator $W_{/ \St}$. 
Later, in \cite{AntCorSto06}, the notion of shorted operator was generalized to that of {\it bilateral shorted operator},  for an operator $W \in L(\HH)$ (not necessarly positive) and a pair of closed subspaces.  The proposed definition comes from the notion of \emph{weak complementability}, which is a refinement of a finite dimensional notion due to T. Ando \cite{Ando}.  In that paper, it was proven that given an operator $W \in L(\HH)$ positive and a closed subspace $\St,$ the shorted operator $W_{/ \St }$ was the maximum of the set $$\M^{-}(W,\St)=\{ X \in L(\HH) : X\minus W, \ R(X) \subseteq \St^{\perp},  \ R(X^*) \subseteq \St^{\perp} \}$$ with the minus order when the pair $(W,\St)$ is compatible.  One of the goals of this paper is to prove that the shorted operator is the maximum of the set $\M^{-}(W,\St)$ with the minus order if and only if the operator $W$ is compatible with the nullspace of $W_{/ \St },$ $N(W_{/ \St }),$ which generalizes the results given in \cite{AntCorSto06}.

Given $W \in L(\HH)^+$ and $A, B \in L(\HH)$ with closed ranges such that $N(B) = N(A^*W)$, in \cite{Contino2} it was shown that, considering the order induced by the cone of positive operators, 
$$
W_{/ R(A)}=\underset{X \in L(\HH)}{inf_{\leq}} \ \ (AXB-I)^*W(AXB-I).
$$
Moreover, also in \cite{Contino2} it was proved the the minimum of the previous set exists if and only if the pair $(W,R(A))$ is compatible. 

In this work we study a similar problem considering the minus order: given $W \in L(\HH)$ a positive operator and $A, B \in L(\HH)$ with closed ranges such that $N(B) = N(A^*W)$, we study the existence of 
\begin{equation}  \underset{X \in L(\HH)}{min_{\minus}} \ \ (AXB-I)^*W(AXB-I),\label{mainproblem} \end{equation}

and its connection with the shorted operator. 

The paper is organized as follows. In section 2, some characterizations of the compatibility of the pair $(W,\St)$ are given. Also the concept of $W$-inverse of an operator $A,$ and some properties are presented. 

In section 3, we collect some useful known results about the minus order, and the connection between compatibility and shorted operators is stated. Also, in this section we prove that the shorted operator of $W$ to $\St,$ denoted by $W_{/ \St},$ is the maximum (in the minus order) of the set $\M^{-}(W,\St)$ if and only if the pair $(W, N(W_{ / \St}))$ is compatible.

In section 4, we study problem \eqref{mainproblem}. We prove that $W_{/R(A)}$ is a lower bound for the set $\{(AXB-I)^{*}W(AXB-I): X \in L(\HH)\}$ (in the minus order) if and only if the pair $(W, N(W_{ / \St}))$ is compatible and in this case, under additional hypothesis, in particular, if $W$ is injective, then $W_{/R(A)}$ is the infimum of that set.

We also prove that $W_{/R(A)}$  is the minimum of the previous set if and only if the pair $(W, R(A))$ is compatible. Moreover, a characterization of the set of solutions is given in terms of the $W$-inverses of $A.$

Motivated by the concepts of the left and the right star orders (see \cite{[BakMit91]}), in section 5 we define the left and right weigthed star orders on $L(\HH).$  These are equivalent to the minus order, with an additional condition on the angle between the ranges of the operators involved imposed by some positive weight. The last part of the section is devoted to appplications. We apply the new characterization of the weighted star order to systems of equations and least squares problems. We also give some formulas for the compression  $W_{\St}= W- W_{/ \St}$ of $W$ to the range of the sum of two closed range operators.

\section{Preliminaries}

In the following $\HH$ denotes a separable complex Hilbert space, $L(\HH)$ is the algebra of bounded linear operators from $\HH$ to $\HH$, and $L(\HH)^{+}$ the cone of semidefinite positive operators. $GL(\HH)$ is the group of invertible operators in $L(\HH),$ $CR(\HH)$ is the subset of $L(\HH)$ of all operators with closed range.
For any $A \in L(\HH),$ its range and nullspace are denoted by $R(A)$ and $N(A)$, respectively. Finally,  $A^{\dagger}$ denotes the  Moore-Penrose inverse of the operator $A \in L(\HH).$ 

Given two closed subspaces $\M$ and $\N$ of $\HH,$ $\M \dot{+} \N$ denotes the direct sum of $\M$ and $\N$ and $\M \ominus \N = \M \cap (\M\cap\N)^{\perp}$. 
If $\HH$ is decomposed as a direct sum of closed subspaces $\HH=\M \dot{+} \N,$ the projection onto $\M$ with nullspace $\N$ is denoted by $P_{\M {\mathbin{\!/\mkern-3mu/\!}} \N},$ and the orthogonal projection onto $\M$ is denoted $P_{\M} = P_{\M {\mathbin{\!/\mkern-3mu/\!}} \M^{\perp}}.$ Also, $\Q$ denotes the subset of $L(\HH)$ of oblique projections, i.e. $\Q=\{Q \in L(\HH): Q^{2}=Q\}.$ 

Given $W\in L(\HH)^{+}$ and a closed subspace $\St$ of $\HH,$ the pair $(W,\St)$ is $\it{compatible}$ if there exists $Q\in \Q$ with $R(Q)=\St$ such that $WQ=Q^{*}W.$ The last condition means that $Q$ is $W$-Hermitian, in the sense that $\PIW{Qx}{y}=\PIW{x}{Qy},$ for every $x, y \in \HH,$ where $\PIW{x}{y}=\PI{Wx}{y}$ defines a semi-inner product on $\HH.$ 
Thus, if $$\mc{P}(W,\St)=\{ Q \in L(\HH): Q^2=Q, \ R(Q)=\St, \ WQ=Q^*W \},$$ then $(W,\St)$ is compatible if and only if $\mc{P}(W,\St)$ is not empty.

The $W$-orthogonal companion of $\St$ is $$\St^{\perp_{W}}=\{x \in \HH: \PI{Wx}{y}=0, \ y\in \St\}.$$
We will use that $\St^{\perp_{W}}=(W\St)^{\perp}=W^{-1}(\St^{\perp}).$

The concept of compatibility between a positive operator $W\in L(\HH)$ and a closed subspace $\St$, has proved to be useful in several applications such as approximation theory, signal processing, among others, see for instance \cite{CG2011}-\cite{Spline}.  

As it was proved in \cite[Prop.~3.3]{CMSSzeged}, the compatibility of the pair $(W,\St)$ is equivalent to a decomposition of the space in terms of the subspace $\St$ and its $W$-orthogonal companion. This is stated in the following theorem. 

\begin{thm} \label{Comp 1} Given $W\in L(\HH)^{+}$ and a closed subspace $\St \subseteq \HH,$ the following conditions are equivalent:
\begin{enumerate} 
\item  [i)] The pair $(W,\St)$ is compatible, 
\item  [ii)] $\HH=\St+\St^{\perp_{W}}$,
\item [iii)] $(W,\St+N(W))$ is compatible,
\item [iv)] $\HH=\St\dot{+} \ (\St^{\perp_{W}}\ominus\St).$
\end{enumerate}
\end{thm}

\medskip 


Suppose that the pair $(W, \St)$ is compatible and let $\N=\St\cap\St^{\perp_{W}},$ observe that $\N=\St \cap N(W);$  we define the projection, $$P_{W, \St}=P_{\St  {\mathbin{\!/\mkern-3mu/\!}} \St^{\perp_{W}} \ominus \N }.$$

\vspace{0,3cm}

The following theorem, proved in \cite[Theorem~3.1]{Arias}, establishes conditions for the existence of solutions of the equation $AXB=C.$ 

\begin{thm} \label{Teosol} Let $A, B, C \in L(\HH).$ If $R(A), R(B)$ or $R(C)$ are closed, then the equation $AXB=C$ admits a bounded solution if and only if $R(C) \subseteq R(A)$ and $R(C^*) \subseteq R(B^*).$ 
\end{thm}

 In \cite{Mitra} S. K. Mitra and C. R. Rao introduced the notion of  $W$-inverse of a matrix, following these ideas, in \cite{Contino} this notion has been extended to linear operators as follows. 
 
 \begin{Def} Given $A \in CR(\HH),$ $B \in L(\HH)$ and $W\in L(\HH)^{+},$ $X_0 \in L(\HH)$ is a $W$-inverse of $A$ in $R(B),$ if for each $x \in \HH$, $X_0x$ is a {\it weighted least squares solution} of $Az=Bx,$ i.e. $$\NWW{AX_0x-Bx}\leq\NWW{Az-Bx}, \mbox{ for every } x, z \in \HH.$$ \end{Def} 
 
 When $B=I,$ $X_0$ is called a $W$-inverse of $A$, see \cite{WGI}. The next theorem characterizes the $W$-inverses.
 
 \begin{thm}
 	\label{thmWinversa} Given $A\in CR(\HH), B \in L(\HH)$ and $W\in L(\HH)^{+},$ the following conditions are equivalent:
 	\begin{enumerate} 
 		\item [i)] The operator $A$ admits a $W$-inverse in $R(B),$
 		\item [ii)] $R(B) \subseteq R(A) + R(A)^{\perp_{W}},$ 
 		\item [iii)] the normal equation $A^{*}W(AX-B)=0$ admits a solution. 
 	\end{enumerate}
 In particular, the pair $(W,R(A))$ is compatible if and only if the normal equation $A^*W(AX-I)=0$ admits a solution.
 \end{thm}
 
 \begin{Dem}
See \cite[Theorem 2.4]{Contino}
 \end{Dem}
 
  \begin{thm}
 	\label{propWinversa} Given $A\in CR(\HH)$ and $W\in L(\HH)^{+}.$ If the pair $(W,R(A))$ is compatible then, for $X_0 \in L(\HH),$ the following conditions are equivalent:
 	\begin{enumerate} 
 		\item [i)] $X_0$ is a $W$-inverse of $A,$
 		\item [ii)] $(AX_0-I)^*W(AX_0-I)=\underset{X \in L(\HH)}{min_{\leq}} \ \ (AX-I)^*W(AX-I)=W_{/ R(A)},$
 		\item [iii)] $X_0$ is a solution of the normal equation $A^{*}W(AX-I)=0.$ 
 	\end{enumerate}
 \end{thm}
 
 \begin{Dem}
 	See \cite[Proposition 4.4]{Contino}.
 \end{Dem}

Section \ref{shorted-minus} is devoted to establish a link between the minus order and the shorted operator. For this purpose, the concept of compatibility will be used.
 
 \section{Shorted operator, compatibility and minus order}\label{shorted-minus}
  
Given a positive operator $W\in L(\HH)^{+}$ and a closed subspace $\St \subseteq \HH$ the notion of shorted operator of $W$ to $\St,$ was introduced by M. G. Krein in \cite{Krein} and later rediscovered by W. N. Anderson and G. E. Trapp, who established a new characterization for this operator  \cite{Shorted1}. As it was proved in \cite{Shorted2}, the set
 $$\M(W,\St)=\{ X \in L(\HH): \ 0\leq X\leq W \mbox{ and } R(X)\subseteq \St^{\perp}\}$$ has a maximum element, with respect to the order $\leq,$ induced in $L(\HH)$ by the cone of positive operators. Then, the {\it shorted operator} of $W$ to $\St$ is defined by
 	$$W_{/\St}=max_{\leq} \ \M(W,\St).$$
 
 	The $\St${\it -compression} $W_{\St}$ of $W$ is the (positive) operator defined by $$W_{\St}=W-W_{/\St}.$$

 For many results on shorted operators, the reader is referred to \cite{Shorted1} and \cite{Shorted2}. Next we collect some results regarding $W_{/\St}$ and $W_{\St}$ which will be useful in the rest of this work. If $W \in L(\HH)^{+},$ $W^{1/2}$ denotes the (positive) square root of $W.$

\begin{thm} \label{TeoShorted}
	Let $W\in L(\HH)^{+}$ and $\St \subseteq \HH$ a closed subspace. Then  
	\begin{enumerate} 
		\item [i)] $W_{/\St}=\mbox{ inf}_{\leq} \ \{ E^{*}WE: E \in \Q, \ N(E)=\St\},$
		\item [ii)] $W_{/\St}=W^{1/2}P_{W^{-1/2}({\St}^{\perp})}W^{1/2},$
		\item [iii)] $R(W) \cap \St^{\perp} \subseteq R(W_{/\St}) \subseteq R(W^{1/2}) \cap \St^{\perp}$, and  $N(W_{/\St})=W^{-1/2}(\overline{W^{1/2}(\St)}),$
		\item [iv)]  $W(\St) \subseteq R(W_{\St}) \subseteq \overline{W(\St)},$ and $N(W_{\St})=W^{-1}(\St^{\perp}).$
	\end{enumerate}

\end{thm}
The formula in $ii)$ was stated by Pekarev, see \cite{Pekarev}.
The inclusions in $iii)$ and $iv)$ can be strict, see \cite{Shorted3}.

Since the infimum appearing in $i)$ of Theorem \ref{TeoShorted} is not attained, it is useful to establish a condition when it is. This is given in the following theorem, see \cite{CMSSzeged} and \cite{Shorted3}.

\begin{thm} \label{TeoShorted2} Let $W\in L(\HH)^{+}$ and $\St \subseteq \HH$ be a closed subspace. The following conditions are equivalent:
	\begin{enumerate} 
		\item [i)] The pair $(W,\St)$ is compatible,
		\item [ii)] $W_{/\St}=\mbox{min } \{ E^{*}WE: E \in \Q, \ N(E)=\St\},$
		\item [iii)] $R(W_{/\St})=R(W) \cap \St^{\perp} \mbox{ and } N(W_{/\St})=N(W)+\St.$
	\end{enumerate}

If any (and then all) of the above conditions holds, then $$W_{/\St}=W(I-Q), \mbox{ for any } Q \in \mc{P}(W,\St).$$
\end{thm}

 \subsection*{\textbf{Minus order and compatibility}}

Different (but equivalent) definitions where given for minus order, for example, using generalized inverses in the matrix case, see \cite{Mit86}. For operators $A, B \in L(\HH),$ in \cite{AntCorSto06}, $A \minus B$ if the \emph{Dixmier} angle between $\ol{R(A)}$ and $\ol{R(B-A)}$ and the Dixmier angle between $\ol{R(A^*)}$ and $\ol{R(B^*-A^*)}$ are less than 1. In this work we give the following definition, equivalent to those appearing in \cite{Sem10} and \cite{Minus}, where in \cite[Theorem 3.3]{Minus},  a characterization of the minus order in terms of the range additivity property is given.


\begin{Def}
	Consider $A, B\in L(\HH)$, we write $A \minus B$  if there exist  (oblique) projections $P, Q \in \Q$  such that $A=PB$ and $A^*=QB^*$.
\end{Def}
	
The projections $P$ and $Q$ can be taken such that $R(P)=\ol{R(A)}$ and $R(Q)=\ol{R(A^*)}.$
It was proven in \cite{Sem10} and \cite{Minus} that $\minus$ is a partial order, known as the minus order for operators.

\begin{thm}\label{equiv minus and RAP}
Consider $A, B\in L(\HH)$. Then the following assertions are equivalent:
	\begin{enumerate}
		\item [i)] $A\minus B,$
		\item [ii)] $R(B)=R(A)\dot{+}R(B-A)$ and $R(B^*)=R(A^*)\dot{+}R(B^*-A^*)$.
	\end{enumerate}
\end{thm}

\bigskip
The left minus order was defined in \cite{Minus} for operators in $L(\HH).$ Using Theorem \ref{equiv minus and RAP}, it is easy to see that this notion is weaker than the minus order, in the infinite dimensional setting. 

\begin{Def}
	Consider $A, B\in L(\HH)$, we write  $A \lminus B$  if $R(B)=R(A)\dot{+}R(B-A)$.
\end{Def}

If $R(B)$ is closed, the left minus order and the minus order are equivalent, in the sense that $A \lminus B$ if and only if $A \minus B,$ as shows the following result, see \cite[Theorem 3.14]{Minus}.

\begin{thm}\label{leftminus then minus}
	Let $A, B\in L(\HH)$ such that $A \lminus B$. If $R(B)$ is closed then $R(A)$ and $R(B-A)$ are closed and $A\minus B$.
\end{thm}

Recall that  $\M(W,\St)=\{ X \in L(\HH): \ 0\leq X\leq W \mbox{ and } R(X)\subseteq \St^{\perp}\}.$

\begin{lema} \label{lemaShorted} Let $W \in L(\HH)^+$ and $\St \subseteq \HH$ be a closed subspace. Then
	\begin{enumerate}
		\item [i)] $\mc{M}(W,N(W_{/ \St})) = \M(W,\St),$
		\item [ii)] $W_{ / \St}= W_{ / N(W_{ / \St}) },$
		\item [iii)] $W(N(W_{ / \St}))^{\perp}=W(\St)^{\perp}.$
	\end{enumerate}
\end{lema}

\begin{dem} $i):$ The inclusion $\St \subseteq N(W_{/ \St})$ (see $iii)$ of Theorem \ref{TeoShorted}) implies that $\mc{M}(W,N(W_{ / \St})) \subseteq \M(W,\St).$ To see the opposite inclusion, take $X \in \mc{M}(W, \St),$ then $ 0 \leq X \leq W_{ / \St}.$ In this case, applying Douglas' theorem \cite{Douglas}, 
$R(X^{1/2}) \subseteq R(W_{ / \St}^{1/2}) \subseteq \ol{R(W_{/ \St})} = N(W_{ / \St})^{\perp}.$ But, $R(X) \subseteq R(X^{1/2}) $ and then $X \in \mc{M}(W,N(W_{ / \St})).$

$ii):$ Using item $i),$ we have that 
$$W_{ / N(W_{/ \St}) } = max_{\leq} \ \mc{M}(W,N(W_{ / \St}))  = max_{\leq } \ \mc{M}(W, \St ) = W_{ / \St}.$$

$iii):$ It follows from item $ii)$ that $W_{  \St}= W_{ N(W_{/ \St}) }$ so that, from $iv)$ of Theorem \ref{TeoShorted}, 
$$W(N(W_{/ \St}))^{\perp}=N(W_{ N(W_{/ \St}) })=N(W_{ \St})=W(\St)^{\perp}.$$

\end{dem}

The next result shows that  the inequality $W_{/ \St} \minus W$ is equivalent to the range inclusion $R(W_{ / \St}) \subseteq R(W)$ and also to a compatibility condition. 

\vspace{0,3cm}

\begin{prop} \label{PropT} Let $W \in L(\HH)^+$ and $\St \subseteq \HH$ be a closed subspace. 
Then the following conditions are equivalent:
\begin{enumerate}
	\item [i)] $W_{/ \St } \minus W,$
	\item [ii)] the pair $(W, N(W_{/ \St}))$ is compatible,
	\item [iii)] $\HH=N(W_{/ \St}) + N(W_{ \St}),$
	\item [iv)] $R(W_{/ \St}) \subseteq R(W).$
 \end{enumerate}

\end{prop}


\begin{dem}
$i) \Leftrightarrow ii):$ Suppose that $W_{/ \St} \minus W,$ then there exists $E \in \Q$ with $R(E^*)=\ol{R(W_{/ \St})},$ such that
$$W_{/ \St}=E^*W=WE.$$ Then $(I-E^*)W=W(I-E)$ and $R(I-E)=N(W_{/ \St})$ Hence, the pair $(W, \TT)$ is compatible.

Conversely, suppose that the pair $(W, \TT)$ is compatible. Let $Q \in \mc{P}(W,\TT)$ then by Theorem \ref{TeoShorted2} and Lemma \ref{lemaShorted},
$$W(I-Q)=\WN=W_{/ \St},$$ and $W_{/ \St} \minus W.$

$ii) \Leftrightarrow iii):$ By Theorem \ref{Comp 1}, $(W, \TT)$ is compatible if and only if
$$\HH=\TT+ \TT^{\perp_{W}}=\TT+ \St^{\perp_{W}}=N(W_{/ \St}) + N(W_{ \St}),$$
where we used Lemma \ref{lemaShorted} and Theorem \ref{TeoShorted}.

$i) \Leftrightarrow iv):$  Suppose that $W_{/ \St } \minus W$, applying Theorem \ref{equiv minus and RAP}, we have that $R(W_{/ \St}) \subseteq R(W).$

Conversely, suppose that $R(W_{/ \St}) \subseteq R(W),$ then $R(W)=R(W_{/ \St})+ R(W_{\St}).$

To see that $R(W_{/ \St}) \cap R(W_{\St})=\{0\},$ consider $x \in  R(W_{/ \St}) \cap R(W_{\St}).$ Since $R(W_{\St}) \subseteq R(W),$ it follows from $iii)$ of Theorem \ref{TeoShorted} that $R(W_{/ \St})=W \cap \St^{\perp}.$ Then there exists $y$ such that $x=Wy$ and $x \in \St^{\perp}.$ Also, from $iv)$ of Theorem \ref{TeoShorted}, $x \in \ol{W(\St)}$ so that there exists a sequence $\{s_n\}_{n \in \mathbb{N}} \subseteq \St $ such that $Ws_n \underset{n \ra \infty }{\ra} x.$ In this case $0=\PI{x}{s_n}=\PI{y}{Ws_n} \underset{n \ra \infty }{\ra} \PI{y}{x}.$ Hence, $0=\PI{y}{x}=\Vert W^{1/2} y \Vert^2,$ and then $W^{1/2}y=0$ and $x=Wy=0.$ Therefore $R(W) = R(W_{/ \St}) \dotplus R(W_{\St}),$ and applying Theorem \ref{equiv minus and RAP}, it follows that $W_{/ \St} \minus W.$

\end{dem}

\begin{cor} \label{CorT} Let $W \in L(\HH)^+$ and $\St \subseteq \HH$ be a closed subspace.  Then the pair $(W, \St)$ is compatible if and only if $W_{/ \St } \minus W$ and $W^{1/2}(\St)$ is closed in $R(W^{1/2}).$
\end{cor}

\begin{dem} If $(W,\St)$ is compatible then, by Theorem \ref{TeoShorted2}, $\TT=\St + N(W).$ By Theorem \ref{Comp 1}, $(W,\TT)$ is also compatible. Applying Proposition \ref{PropT}, it follows that $W_{/ \St} \minus W.$ 

Also, $(W,\St)$ is compatible if and only if $\HH=\St + W^{-1} (\St^{\perp})$ and applying $W^{1/2}$ to both sides of the equality, we get 
$$R(W^{1/2})=W^{1/2} \St \oplus (W^{1/2} \St)^{\perp} \cap R(W^{1/2}) \subseteq \ol{W^{1/2} \St }	\cap R(W^{1/2}) \oplus (W^{1/2} \St)^{\perp} \cap R(W^{1/2}) \subseteq R(W^{1/2}).$$
Therefore, $W^{1/2} \St =  \ol{W^{1/2} \St }	\cap R(W^{1/2}),$ so that $W^{1/2}(\St)$ is closed in $R(W^{1/2}).$ See also, \cite[Proposition 3.8]{Shorted3}.

Conversely, suppose that $W_{/ \St} \minus W$ and $W^{1/2}(\St)$ is closed in $R(W^{1/2}).$ Then, by Proposition \ref{PropT}, the pair $(W,\TT)$ is compatible. But, by Theorem \ref{TeoShorted},
$$\TT=W^{-1/2}(\ol {W^{1/2}(\St)} \cap R(W^{1/2}))=W^{-1/2}(W^{1/2}(\St))=\St+N(W).$$ 
Therefore, by Theorem \ref{Comp 1}, the pair $(W, \St)$ is compatible.
\end{dem}

\vspace{0,3cm}

The shorted operator can also be  characterized  as the maximum of certain set when the minus order is considered. More precisely, for $W \in L(\HH)^+$ and $\St \subseteq \HH$ a closed subspace, define
$$\M^{-}(W,\St)=\{ X \in L(\HH) : X\minus W, \ R(X) \subseteq \St^{\perp},  \ R(X^*) \subseteq \St^{\perp} \}.$$

In \cite{Mit86}, Mitra proved (for matrices in $\mathbb{C}^{m \times n}$) that the shorted operator is the maximum of the set $\M^{-}(W,\St),$ where the partial ordering is the minus order. In \cite{AntCorSto06}, Antezana et al. proved a similar result for operators when the pair $(W,\St)$ is compatible. The next results generalize this fact.

\begin{lema} \label{lemaMaximo} Let $W \in L(\HH)^+$ and $\St \subseteq \HH$ be a closed subspace. Then $$\M^{-}(W,\TT)=\M^{-}(W,\St).$$
\end{lema}

\begin{dem} Since $\St \subseteq \TT,$ then $\TT^{\perp} \subseteq \St^{\perp}$ and  $\M^{-}(W,\TT)\subseteq \M^{-}(W,\St).$
	
	On the other hand, let $X \in  \M^{-}(W,\St),$ then $X \minus W$ and so, by Theorem \ref{equiv minus and RAP}, $R(X) \subseteq R(W).$
	Therefore, if $R(X) \subseteq \St^{\perp}$ then $W^{-1} (R(X)) \subseteq W^{-1} (\St^{\perp})= W^{-1}(\TT^{\perp})$. Then
	$$W(W^{-1} (R(X)))  \subseteq \TT^{\perp},$$ but $W(W^{-1}(R(X))) = R(X) \cap R(W) = R(X).$ Hence, $R(X) \subseteq \TT^{\perp}.$
	Analogously, $R(X^*) \subseteq \TT^{\perp},$ because, since $W$ is positive, $X \minus W$ implies that $X^* \minus W$.
	Therefore $\M^{-}(W,\St)\subseteq \M^{-}(W,\TT).$
	
\end{dem}

\begin{thm} \label{thmMaximo} Let $W \in L(\HH)^+$ and $\St \subseteq \HH$ be a closed subspace. Then the pair $(W, N(W_{ / \St}))$ is compatible if and only if $$W_{/ \St } = max_{\minus} \ \ \M^{-}(W,\St).$$ 
\end{thm}

\begin{dem} Suppose the pair $(W, N(W_{ / \St}))$ is compatible. By Proposition \ref{PropT}, $W_{/ \St} \minus W.$ By Lemma \ref{lemaShorted},  $W_{/ \St} = W_{/ \TT},$ therefore $R(W_{/ \St})=R(W_{/ \TT}) \subseteq \TT^{\perp}.$ Hence  $W_{/ \St} \in \M^{-}(W,\TT).$
	
	On the other hand, given $X \in \M^{-}(W,\TT),$ from $X \minus W,$ there exists $E \in \Q$ such that $X=EW.$ Let $Q \in \mc{P}(W, \TT)$ then by Theorem \ref{TeoShorted2}, $W_{/ \TT}=W(I-Q).$
	The inclusion $R(X^*) \subseteq \TT^{\perp} = R(Q)^{\perp}=N(Q^*),$ implies that $Q^*X^*=0,$ or equivalently $X(I-Q)=X.$ Then
	
	$$X=X(I-Q)=EW(I-Q)=EW_{/ \TT}=EW_{/ \St}.$$ In a similar way, there exists a projection $F$ such that $X^*=FW_{/ \St}.$ Therefore $X \minus W_{/ \St}.$
	
	Then, by Lemma \ref{lemaMaximo}, $$W_{/ \St } = max_{\minus} \ \ \M^{-}(W,\TT)=max_{\minus} \ \ \M^{-}(W,\St).$$   
	
	Conversely, if $W_{/ \St } =max_{\minus} \ \ \M^{-}(W,\St),$ in particular $ W_{/ \St } \minus W,$ and by Proposition \ref{PropT}, the pair $(W,\TT)$ is compatible.
	
\end{dem}

\section{Shorted operator characterizations and the minus order}

Let $W \in L(\HH)^+$ and $A, B\in CR(\HH)$ such that $N(B) =  R(A)^{\perp_{W}}.$ 

In \cite[Proposition 3.1]{Contino2} it was proved that the infimum (in the order induced by the cone of positive operators) of the set $\{ (AXB-I)^*W(AXB-I) :  X \in L(\HH) \}$ exists and   
\begin{equation}  \label{InfimoOP1}
\underset{X \in L(\HH)}{inf_{\leq}} \  (AXB-I)^*W(AXB-I) = W_{/ R(A)}.
\end{equation}
Also in \cite[Theorem 3.2]{Contino2}, it was proved that the minimum of this set exists if and only if the pair $(W,R(A))$ is compatible. 

In this section, we study a similar problem considering the minus order: for $W \in L(\HH)^+$ and $A, B\in CR(\HH),$ analyze the existence of 
%
%

\begin{equation}  \label{EqProblem}
\underset{X \in L(\HH)}{min_{\minus}} \  (AXB-I)^*W(AXB-I).
\end{equation}

Problem \eqref{EqProblem} can be restated as a minimization problem with a constraint. In fact, problem \eqref{EqProblem} is equivalent to the following problem: analize the existence of

$$\underset{X \in L(\HH)}{min_{\minus}} \  (AX-I)^*W(AX-I) \mbox{ subject to } \ N(A^*W) \subseteq N(X).$$

\bigskip

\begin{prop} \label{Propminus2} Let $A, B \in CR(\HH)$ and $W \in L(\HH)^+$  such that $N(B)=R(A)^{\perp_{W}}.$ Then the shorted operator $W_{/ R(A)}$ is a lower bound for $\{ (AXB-I)^*W(AXB-I) : \ X \in L(\HH)\}$ if and only if the pair $(W, N(W_{/ R(A)}))$ is compatible.

In this case, if $R(A)+ R(A)^{\perp_{W}}$ is closed, then 
$$W_{ / R(A)}= \underset{X \in L(\HH)}{inf_{\minus}} \  (AXB-I)^*W(AXB-I).$$

\end{prop}

\begin{dem}
	Let $X \in L(\HH),$ then
	$$F(X): =(AXB-I)^*W(AXB-I)=W_{/ R(A)} + (AXB-I)^*W_{R(A)} (AXB-I).$$
	By Lemma \ref{lemaShorted}, we have that $W_{ / \TTA}=W_{/ R(A)}$ and $\TTA^{\perp_{W}}=R(A)^{\perp_{W}}.$
	
	Suppose  the pair $(W, \TTA)$ is compatible, let $Q \in \mc{P}(W, \TTA),$
	then
	$$F(X)(I-Q)=W_{/ R(A)}(I-Q)+(AXB-I)^*W_{R(A)}(AXB-I)(I-Q)=W_{/ R(A)},$$
	where we used that $R(Q)=N(W_{/ R(A)}),$ then $W_{/R(A)}(I-Q)=W_{/R(A)},$ and the facts that 
	$R(I-Q)=N(Q)\subseteq \TTA^{\perp_{W}} = R(A)^{\perp_{W}}= N(B)$ and $N(W_{R(A)}) = R(A)^{\perp_{W}}.$ Therefore
	$$F(X)(I-Q)=W_{/ R(A)}=(I-Q^*)F(X), \mbox{ for every } X \in L(\HH).$$  Then
	$$W_{/ R(A)} \minus F(X), \mbox{ for every } X \in L(\HH).$$
	Hence, $W_{/ R(A)}$ is a lower bound for $\{ F(X) : \ X \in L(\HH)\}.$
	
	Conversely, suppose $W_{/ R(A)}$ is a lower bound for $\{ F(X) : \ X \in L(\HH)\},$ in particular $W_{/ R(A)} \minus W,$ and by Proposition \ref{PropT}, the pair $(W, \TTA)$ is compatible.
	
	Finally, suppose that $(W, \TTA)$ is compatible and  $R(A)+ R(A)^{\perp_{W}}$ is closed. 
	Then $R(A)+ N(B)$ is closed and by \cite[Corollary 2.5]{Izumino} $R(BA)$ is closed. 
	Also, since $N(B)=N(A^*W),$ we have that $N(BA)=N(A^*WA).$ Therefore  $R(A^{*}W^{1/2}) \subseteq \ol{R(A^{*}WA)}=\ol{R(A^*B^*)}= R(A^*B^*).$ Then, by Theorem \ref{Teosol}, there exists $X_0 \in L(\HH)$ such that
	$W^{1/2}AX_0BA=W^{1/2}A.$ Then $R(A) \subseteq N(W^{1/2}(AX_0B-I))=N(F(X_0)),$ or equivalently 
	$$R(F(X_0)) \subseteq R(A)^{\perp}.$$
	Let $D\in L(\HH)$ be any lower bound for $\{F(X) : X \in L(\HH) \},$ then
	$$D \minus F(X), \mbox{ for every } X \in L(\HH).$$ In particular, $D \minus W.$ 
	 Since $D\minus F(X_0),$ by Theorem \ref{equiv minus and RAP}, $R(D) \subseteq R(F(X_0)) \subseteq R(A)^{\perp}.$ From $D \minus F(X_0)$ and the fact that $F(X_0)$ is positive we get that $D^* \minus F(X_0).$ Therefore, in the same way as before, we get 
	$R(D^*) \subseteq R(A)^{\perp}.$ Then $D \in \M^{-}(W,\St)$ and since $(W, \TTA)$ is compatible, by Theorem \ref{thmMaximo}, $D \minus W_{/ R(A)}.$
\end{dem}

\bigskip
For example, if $W$ is injective, then 
$W_{ / R(A)}= \underset{X \in L(\HH)}{inf_{\minus}} \  (AXB-I)^*W(AXB-I).$
\bigskip
In fact, if $W$ is injective and the pair $(W, N(W_{/ R(A)}))$ is compatible, by Theorem \ref{Comp 1}, $\HH=\TTA \dotplus W^{-1}(\TTA^{\perp})=\TTA \dotplus W^{-1}(R(A)^{\perp}),$ and since $R(A) \subseteq \TTA,$ it follows that $R(A) \dotplus W^{-1}(R(A)^{\perp})$ is also closed and by the last theorem we get the result.

\bigskip
The next proposition shows that the minimum of \eqref{EqProblem} in the minus order  is $W_{/ R(A)}$ if and only if the pair $(W, R(A))$ is compatible. First, we need the following lemma which shows that when $W_{/R(A)}$ is in the image of the function $G(X)=(AX-I)^*W(AX-I),$ the pair $(W,R(A))$ is compatible. For its proof, we follow the same ideas as in \cite[Theorem 3.2]{Contino}.

\begin{lema} \label{LemaShorted} Let $A \in CR(\HH)$ and $W \in L(\HH)^+,$  there exists $X_0 \in L(\HH)$ such that $$W_{/R(A)}=(AX_0-I)^*W(AX_0-I),$$if and only if the pair $(W,R(A))$ is compatible.
\end{lema}

\begin{dem}
	Suppose that $W_{/R(A)}=(AX_0-I)^*W(AX_0-I),$ for certain $X_0 \in L(\HH).$ Writing $W=W_{/R(A)}+W_{R(A)}$, it follows that
	$$W_{/R(A)}=(AX_0-I)^*W(AX_0-I)=W_{/R(A)}+(AX_0-I)^{*}W_{R(A)}(AX_0-I),$$  because $R(A)\subseteq N(W_{/ R(A)}).$
	Therefore $$(AX_0-I)^{*}W_{R(A)}(AX_0-I)=0, \mbox{ or, equivalently, } W_{R(A)}^{1/2}(AX_0-I)=0.$$ 
	Then, by $iv)$ of Theorem \ref{TeoShorted},
	$$R(AX_0-I) \subseteq N(W_{R(A)})=W^{-1}(R(A)^{\perp}), \mbox{ or }W(R(AX_0-I))\subseteq R(A)^{\perp} \cap R(W).$$
	Then $$R(W_{/R(A)})=R((AX_0-I)^*W(AX_0-I))\subseteq (AX_0-I)^{*} (R(A)^{\perp} \cap R(W))=R(A)^{\perp} \cap R(W),$$
	because $A^{*}(R(A)^{\perp})=0.$ Then, by $iii)$ of Theorem \ref{TeoShorted}, $R(W_{/R(A)})= R(A)^{\perp} \cap R(W).$ 
	
	Also $x \in N(W_{/R(A)})$ if and only if $W^{1/2}(AX_0-I)x=0,$ or equivalently $(AX_0-I)x \in N(W).$
	In this case $x \in N(W)+R(A),$ and then again applying $iii)$ of Theorem \ref{TeoShorted}, $$N(W_{/R(A)})=N(W)+R(A).$$ 
	Therefore 
	$R(W_{/R(A)}) = R(A)^{\perp} \cap R(W)$ and $N(W_{/R(A)})=N(W)+R(A)$ and by Theorem \ref{TeoShorted2}, the pair $(W, R(A))$ is compatible.
	
	Conversely, if the pair $(W, R(A))$ is compatible, by Theorem \ref{TeoShorted2}, $W_{/R(A)}=W(I-Q)=(I-Q)^*W(I-Q),$ for any $Q \in \mc{P}(W,R(A)).$ Consider $X_0=A^{\dagger}Q,$ then $AX_0=Q$ and $F(X_0)=W_{/R(A)}.$ 
\end{dem}

\bigskip

 
\medskip
We denote $S_p, $ for some $1 \leq p < \infty,$ the p-Schatten class, see \cite{Ringrose}.
Let $W \in L(\HH)^+,$ such that $W^{1/2} \in S_p$ for some $1 \leq p < \infty,$ in the following we define 
$$\Vert X\Vert_{p,W}=\Vert W^{1/2} X \Vert_{p},$$ for any $X \in L(\HH).$ For $A, B \in CR(\HH)$ such that $N(B)=R(A)^{\perp_{W}},$  the next proposition, shows the equivalence between the existence of minimum of problem \eqref{EqProblem} and the existence of
$\underset{X \in L(\HH)}{min} \Vert AXB-I \Vert_{p,W}.$

\begin{prop} \label{Propminus} Let $A, B \in CR(\HH)$ and $W \in L(\HH)^+$  such that $N(B)=R(A)^{\perp_{W}}.$ Then the following are equivalent:

\begin{enumerate}
\item [i)]	 $W_{/ R(A)}=\underset{X \in L(\HH)}{min_{\minus}} \ (AXB-I)^*W(AXB-I),$ 
\item [ii)]  the pair $(W, R(A))$ is compatible,
\item [iii)] $W_{/ R(A)}=\underset{X \in L(\HH)}{min_{\leq}}  \ (AXB-I)^*W(AXB-I).$ 
\end{enumerate}

If $W^{1/2} \in S_p$ for some $1 \leq p < \infty,$ conditions $i), ii)$ and $iii)$ are also equivalent to

\begin{enumerate}
\item [iv)] $\Vert W_{/ R(A)}^{1/2} \Vert_{p} = \underset{X \in L(\HH)}{min} \Vert AXB-I \Vert_{p,W}.$
\end{enumerate}

\end{prop}

\begin{dem}  $i) \Leftrightarrow ii):$ Suppose the pair $(W, R(A))$ is compatible. Then, by Theorem \ref{Comp 1} and Theorem \ref{TeoShorted2}, $\TTA=R(A)+N(W)$ and the pair $(W,\TTA)$ is compatible. Thus, by Proposition \ref{Propminus2}, $W_{/ R(A)}$ is a lower bound for $\{ F(X) : X \in L(\HH) \},$ where $F(X):= (AXB-I)^*W(AXB-I).$ It remains to see that there exists $X_0 \in L(\HH)$ such that $W_{/ R(A)}=F(X_0).$
	
Let $Q \in \mc{P}(W, R(A)),$ then  by Theorem \ref{TeoShorted2}, $W_{/ R(A)}=(I-Q)^*W(I-Q).$
Observe that $R(W^{1/2}Q)=R(W^{1/2}A).$ On the other hand, since $B \in CR(\HH)$ and $N(W^{1/2}Q)=N(WQ)=N(Q^*W)=W^{-1}(N(Q^*))=W^{-1}(R(Q)^{\perp})= R(A)^{\perp_{W}}=N(B),$ we have that $R(Q^*W^{1/2}) \subseteq N(W^{1/2}Q)^{\perp} = N(B)^{\perp} = R(B^*).$ Then, by Theorem \ref{Teosol}, there exists $X_0 \in L(\HH)$, such that $W^{1/2}AX_0B=W^{1/2}Q.$
Thus 
$$(AX_0B-I)^*W(AX_0B-I)=(I-Q)^*W(I-Q) = W_{/ R(A)},$$ and $W_{/ R(A)} \in \{ F(X) : X \in L(\HH) \}.$ Hence, $$W_{/ R(A)}=\underset{X \in L(\HH)}{min_{\minus}} \ F(X).$$
	
Conversely, suppose $W_{/ R(A)}=\underset{X \in L(\HH)}{min_{\minus}} \ F(X).$ Then, there exists $X_0\in L(\HH),$ such that
$$(AX_0B-I)^*W(AX_0B-I)=\underset{X \in L(\HH)}{min_{\minus}} \ F(X) =W_{/R(A)}.$$ Then, by Lemma \ref{LemaShorted}, the pair $(W,R(A))$ is compatible. 

$ii) \Leftrightarrow iii)$ and $ii) \Leftrightarrow iv):$ follow from the fact that $N(B) = N(A^*W)$ and \cite[Theorem 4.3]{Contino2}.
	
\end{dem}

\begin{prop} \label{Propmin} Let $A, B \in CR(\HH)$ and $W \in L(\HH)^+$  such that $N(B)=R(A)^{\perp_{W}}$ and the pair  $(W, R(A))$ is compatible. Then, for $X_0 \in L(\HH),$ the following conditions are equivalent:
	\begin{enumerate}
		\item [i)]  $(AX_0B-I)^*W(AX_0B-I)=\underset{X \in L(\HH)}{min_{\minus}} \ (AXB-I)^*W(AXB-I) =W_{/ R(A)},$
		\item [ii)]  $A^*W(AX_0B-I)=0,$
		\item [iii)] $X_0B$ is a $W$-inverse of $A.$ 
	\end{enumerate}
 
\end{prop}

\begin{dem} The equivalence between $i), \ ii), \ iii)$ follows from the fact that $X_0 \in L(\HH)$ satisfies 
	$(AX_0B-I)^*W(AX_0B-I)=\underset{X \in L(\HH)}{min_{\minus}} \ (AXB-I)^*W(AXB-I)=W_{/ R(A)}$ if and only if $X_0$ satisfies $(AX_0B-I)^*W(AX_0B-I)=\underset{X \in L(\HH)}{min_{\leq}} \ (AXB-I)^*W(AXB-I)=W_{/ R(A)}$ (see Proposition \ref{Propminus}) and \cite[Theorem 3.2]{Contino2}.
\end{dem}

\section{Weighted star order and applications}

In  \cite{[Dra78]}, Drazin introduced  the  star order on semigroups with involutions and in \cite{[BakMit91]}, Baksalary and Mitra defined the left and right star orders on the set of complex matrices. Later, Antezana et al. studied the star  order  on the algebra of bounded operators on a Hilbert space \cite{[AntCanMosSto10]}.

Given $A, B\in L(\HH),$  the {\it star order} $A \star B$, the {\it left star order} $A \lstar B$ and  the {\it right star order} $A \rstar B$ are defined, respectively, by  $A\star B$ if 
$
A^*A=A^*B \textrm{ and } AA^*=BA^*;
$
$A \lstar B$ if 
$
A^*A=A^*B  \textrm{ and } R(A)\subseteq R(B) 
$
and $A \rstar B$  if 
$
AA^*=BA^*  \textrm{ and } R(A^*)\subseteq R(B^*).
$

\medskip
If $A, B\in L(\HH)$, then $A\star B$ if and only if there exist orthogonal projections $P, Q$ such that $A=PB$ and  $A^*=QB^*$, see  \cite[Proposition 2.3]{[AntCanMosSto10]}.

%
%
%
%
\medskip
 Observe that $A\star B$ if and only  $A \minus B,$  $R(B-A) \subseteq R(A)^{\perp}$ and $R(B^*-A^*) \subseteq R(A^*)^\perp.$ More generally, given a positive operator $W,$ we now introduce the weighted star order, which is the minus order with an additional condition on the angle between the ranges of the operators involved, imposed by the weight $W.$

\begin{Def} Given $A,B \in L(\HH)$ and  $W \in  L(\HH)^{+},$ we say $A \lW B$ if 
$A \lminus B$ and $R(B-A)\subseteq R(A)^{\perp_W}.$ 

Analogously, we say $A \rW B$ if $A\rminus B$  and $R(B^*-A^*)\subseteq R(A^*)^{\perp_W}.$ 

Finally, we say $A \starW B$ if $A \lW B$ and $A \rW B.$
\end{Def}
	
\begin{prop} \label{Prop1} Let $A,B \in L(\HH)$ and  $W \in  L(\HH)^{+},$ then $A \lW B$ if and only if $R(B)=R(A)\dotplus R(B-A)$ and $A^*WA=A^*WB.$
\end{prop}

\begin{dem} Suppose that $R(B)=R(A)\dotplus R(B-A)$ and $A^*WA=A^*WB,$ then by the definition of $\lminus$,  $A  \lminus B$ and since $A^*WA=A^*WB$ it follows that $R(B-A) \subseteq N(A^*W)=R(A)^{\perp_W},$ then $A \lW B.$ 
	
Conversely, if  $A \lW B,$ then $A  \lminus B,$ or $R(B)=R(A)\dotplus R(B-A),$ and since $R(B-A) \subseteq R(A)^{\perp_W}=N(A^*W),$ we have $A^*WA=A^*WB.$ 
\end{dem}

\begin{prop} \label{Prop3} Let $W \in  L(\HH)^{+},$ then the relations $\lW,$ $\rW$ and $\starW$ are partial orders. 
\end{prop}

\begin{dem} Observe that $\lW,$ $\rW$ and $\starW$ are clearly reflexive and also antisymmetric, because $\minus,$ $\rminus$ and $\lminus$ are antisymmetric (see \cite[Proposition 3.11]{Minus}). 
	
Finally, if $A \lW B$ and $B \lW C,$ then $A  \lminus B$ and it follows that $R(A) \subseteq R(B),$ therefore $R(W^{1/2}A) \subseteq R(W^{1/2}B).$ Also $A^*WA=A^*WB,$ so that 
$W^{1/2}A \lstar W^{1/2} B.$  In the same way  $W^{1/2}B \lstar W^{1/2} C,$ then since $\lstar$ is transitive,  $W^{1/2}A \lstar  W^{1/2} C$ and therefore $A^*WA=A^*WC.$ On the other hand, if  $A \lW B$ and $B \lW C,$ then $A \lminus B$ and $B \lminus C,$ and since $\lminus$ is transitive  (see \cite[Proposition 3.11]{Minus}), we also have that $A \lminus C.$  Then, by Proposition \ref{Prop1}, $A \lW C$ and the relation $\lW$ is a partial order. In a similar way, it can be proven that $\starW,$ $\rW$ are transitive.
\end{dem}

\begin{lema} \label{LemmaEstrellaW} Let $A, \ B \in L(\HH)$ and  $W \in  L(\HH)^{+},$ such that $N(W) \cap \ol{R(A)} =\{0\}$ and the pair $(W, \ol{R(A)})$ is compatible. Then the following are equivalent:
\begin{enumerate}
	\item [i)] $A \lW B,$
	\item  [ii)] $R(A) \subseteq R(B)$ and $A=P_{W, \ol{R(A)}} B.$ 
	\item [iii)] $R(A) \subseteq R(B)$ and $A^*WA=A^*WB.$
\end{enumerate}
\end{lema}

\begin{dem}
$i) \Rightarrow ii):$ Suppose $A \lW B,$ then $R(B)=R(A)\dotplus R(B-A)$ and $R(B-A)\subseteq R(A)^{\perp_W}=N(P_{W, \ol{R(A)}}),$ therefore $R(A) \subseteq R(B).$ Hence
$0=P_{W, \ol{R(A)}} (B-A)=P_{W, \ol{R(A)}} B-A=0$ and $A=P_{W, \ol{R(A)}} B.$

$ii) \Rightarrow iii):$ Suppose $R(A) \subseteq R(B)$ and $A=P_{W, \ol{R(A)}} B.$ Then
$R(B-A) \subseteq N(P_{W, \ol{R(A)}} ) \subseteq R(A)^{\perp_W}=N(A^*W),$ therefore $A^*WA=A^*WB.$

$iii) \Rightarrow i):$ Suppose $R(A) \subseteq R(B)$ and $A^*WA=A^*WB.$ The inclusion $R(A) \subseteq R(B)$ is equivalent to $R(B)= R(A) + R(B-A),$ see \cite[Proposition 2.4]{AriCorMae15}.

On the other hand, from $A^*WA=A^*WB,$ it follows that $R(B-A) \subseteq N(A^*W)=R(A)^{\perp_W}.$ Finally, if $y \in R(A) \cap  R(B-A),$ then there exists $x, z \in \HH$ such that $$y=Ax=(B-A)z.$$ Multiplying the last equation by $A^*W,$ we get
$$A^*WAx=A^*W(B-A)z=0,$$ then $WAx=0,$  so that $Ax \in \ol{R(A)} \cap N(W)=\{0\},$ hence $y=Ax=0$ and $R(B)= R(A) \dotplus R(B-A).$ Then $A \lW B.$

\end{dem}

A similar result can be stated for $\rW.$

\begin{cor} \label{CorEstrellaW3} Let $A \in CR(\HH), \ B \in L(\HH)$ and  $W \in  L(\HH)^{+}$ such that $N(W) \cap \ol{R(A)} =\{0\},$ $N(W) \cap \ol{R(A^*)} =\{0\}$ and the pairs $(W,\ol{R(A)})$  and $(W,\ol{R(A^*)})$ are compatible. Then the following are equivalent:
	\begin{enumerate}
		\item [i)]  $A\starW B,$
		\item [ii)]  $A^*=P_{W, \ol{R(A^*)}} B^*$  and $A=P_{W, \ol{R(A)}} B,$ 
		\item [iii)] $AWA^*=BWA^*$ and $A^*WA=A^*WB.$
	\end{enumerate}
\end{cor}
	
\begin{Dem} Straightforwards.
\end{Dem}

\begin{thm} \label{THM1} Let $A,B \in L(\HH)$ such that $R(A+B)$ is closed and $C \in L(\HH).$ Let $W \in  L(\HH)^{+}$ such that $(W,R(A+B))$ is compatible, $W^{1/2} \in S_p$ for some $1\leq p < \infty,$ and $A \lW A+B.$ Then the following statements are equivalent:
\begin{enumerate}
	\item [i)] $X_0 \in L(\HH)$ satisfies
	$$\Vert(A+B)X_0-C\Vert_{p,W}= \underset{X \in L(\HH)}{min}  \Vert(A+B)X-C\Vert_{p,W};$$ 
	\item [ii)]  $X_0 \in L(\HH)$ satisfies
	$$\left \{
	\begin{array}{ll}
	\Vert AX_0-C\Vert_{p,W} = \underset{X \in L(\HH)}{min}  \Vert AX-C\Vert_{p,W}\\
	\Vert BX_0-C\Vert_{p,W} = \underset{X \in L(\HH)}{min}  \Vert BX-C\Vert_{p,W}.
	\end{array} \right. $$
\end{enumerate}
\end{thm}

\begin{dem}
First observe that if $R(A+B)$ is closed and $A \lW A+B,$ by Proposition \ref{Prop1} and Corollary \ref{leftminus then minus}, we have $A^*WB=0$ and $A \minus A+B.$

Now, suppose item $i)$ holds, then by \cite[Theorem 4.5]{Contino}, $X_0$ is a solution of
$$(A^*+B^*)W((A+B)X-C)=0.$$ 
Therefore, $R(W(A+B)X_0-C) \subseteq N(A^*+B^*)=N(A^*)\cap N(B^*),$ because $A \minus A+B.$
Then
$$A^*W((A+B)X_0-C)=B^*W((A+B)X_0)-C)=0,$$ and since $A^*WB=0,$ it holds that $$A^*W(AX_0-C)=B^*W(BX_0-C)=0,$$  
or equivalently by \cite[Theorem 4.5]{Contino}, $X_0 \in L(\HH)$ satisfies
$$\left \{
\begin{array}{ll}
\Vert AX_0-C\Vert_{p,W} =  \underset{X \in L(\HH)}{min}  \Vert AX-C\Vert_{p,W}\\
\Vert BX_0-C\Vert_{p,W} =  \underset{X \in L(\HH)}{min}  \Vert BX-C\Vert_{p,W}.
\end{array} \right. $$

Conversely, suppose item $ii)$ holds, then by \cite[Theorem 4.5]{Contino}, $X_0$ is a solution of the system 
$$\left \{
\begin{array}{ll}
A^*W(AX-C)=0\\
B^*W(BX-C)=0,
\end{array} \right.$$
therefore, since $A^*WB=0,$ we have $(A^*+B^*)W((A+B)X_0-C)=A^*W(AX_0-C)+B^*W(BX_0-C)=0.$ Then $X_0$ is as solution of $(A^*+B^*)W((A+B)X-C)=0$ and by \cite[Theorem 4.5]{Contino}, we have item $i).$
\end{dem}

\begin{cor} \label{CorWinversas} Let $A,B \in L(\HH)$ such that $R(A+B)$ is closed. Let $W \in  L(\HH)^{+}$ such that the pair $(W,R(A+B))$ is compatible and $A \lW A+B.$  Then $X_0$ is a $W$-inverse of $A+B$ in $R(C)$ if and only if $X_0$ is a $W$-inverse of $A$ in $R(C)$ and a $W$-inverse of $B$ in $R(C),$ 
\end{cor}

\begin{dem}
	By Theorem \ref{thmWinversa},  $X_0$ is $W$-inverse of $A+B$ in $R(C)$ if and only if $X_0$ is a solution of $$(A^*+B^*)W((A+B)X-C)=0,$$ if and only if, by the proof of Theorem \ref{THM1}, $X_0$ is a solution of the system 
	$$\left \{
	\begin{array}{ll} 
	A^*W(AX-C)=0,\\
	B^*W(BX-C)=0,
	\end{array} \right.$$
	or equivalently, again by Theorem \ref{thmWinversa}, $X_0$
	is a $W$-inverse of $A$ in $R(C)$ and a $W$-inverse of $B$ in $R(C).$
\end{dem}

\begin{cor}  Let $A,B \in L(\HH)$ with $R(A+B)$ closed. Let $W \in  L(\HH)^{+}$ such that $(W,R(A+B))$ is compatible and $A \lW A+B.$ Then
$$W_{R(A+B)}=W_{ R(A)}+W_{ R(B)}.$$
\end{cor}

\begin{dem} It follows from Theorem \ref{thmWinversa}, for $C=I,$ that $(W,R(A+B))$ is 
compatible if and only if there exists $X_0$ such that 
\begin{equation} \label{eq1Cor5.8}
(A^*+B^*)W((A+B)X_0-I)=0.
\end{equation}

By Corollary \ref{CorWinversas}, this is equivalent to 
\begin{equation} \label{eq2Cor5.8}
\left \{
\begin{array}{ll} 
A^*W(AX_0-I)=0,\\
B^*W(BX_0-I)=0.
\end{array} \right.
\end{equation}

On the other hand, by Theorem \ref{propWinversa}, equation \eqref{eq1Cor5.8} and \eqref{eq2Cor5.8} are equivalent to the equalities
 $$W_{/ R(A+B)}=((A+B)X_0-I)^*W((A+B)X_0-I),$$
 $$W_{/ R(A)}=(AX_0-I)^*W(AX_0-I) \mbox{ and } W_{/ R(B)}=(BX_0-I)^*W(BX_0-I).$$
 But, using the fact that $A \lW A+B,$ it follows that
 $$W_{/ R(A+B)}=((A+B)X_0-I)^*W((A+B)X_0-I)=$$
 $$=(AX_0-I)^*W(AX_0-I)+(BX_0-I)^*W(BX_0-I)-W=W_{ / R(A)}+W_{ / R(B)}-W.$$
 Therefore,
 $$W_{R(A+B)}=W_{ R(A)}+W_{ R(B)}.$$
\end{dem}

\section*{Acknowledgements}
This research was partially supported by CONICET PIP 0168.


\end{document}